\input vanilla.sty
\scaletype{\magstep1}
\scalelinespacing{\magstep1}
\def\bull{\vrule height .9ex width .8ex depth -.1ex}


\pageno=1

\title Operators preserving orthogonality are isometries
\endtitle

\author Alexander Koldobsky
\\ Department 
of Mathematics\\ University of Missouri-Columbia\\ Columbia, MO 65211\\
\endauthor

\vskip1truecm

\subheading{Abstract} Let $E$ be a real Banach space. For $x,y \in E,$ we 
follow 
R.James in saying that $x$ is orthogonal to $y$ if $\|x+\alpha y\|\geq \|x\|$
for every $\alpha \in R$. We prove that every operator from $E$ into itself
preserving orthogonality is an isometry multiplied by a constant.

\vskip1truecm

Let $E$ be a real Banach space. For  $x,y \in E,$ we follow R.James in saying 
that $x$ is orthogonal to $y$  $(x\bot y)$  
if $\|x+\alpha y\|\geq \|x\|$ for every 
$\alpha \in R$. 

It is clear that every isometry $T:E\rightarrow E$ preserves 
orthogonality, i.e.
$x\bot y$ implies $Tx\bot Ty.$ We prove here that 
the converse statement is valid,
namely, every linear operator preserving orthogonality is an 
isometry multiplied
by a constant.

D.Koehler and P.Rosenthal  [4]  have proved that an 
operator is an isometry if 
and only if it preserves any semi-inner product. It is easy to show (see [2])
that orthogonality of vectors with respect to any semi-inner product implies
James' orthogonality. So the result of this paper seems to 
refine that from [4].

We start with some auxiliary facts.

For $x\in E, x\neq 0,$ denote by $S(x)=\{x^{*}\in E^{*}: \|x^{*}\|=1, 
x^{*}(x)=\|x\|\}$
the set of support functionals at the point $x$. It is well-known [1] that,
for every  $x,y \in E, x\neq 0,$

$$\lim_{\alpha\rightarrow 0} (\|x+okstate.\alpha y\|-\|x\|)/\alpha = 
\sup_{x^{*}\in S(x)}\ x^{*}(y) \tag 1$$

The limit in the left-hand side as $\alpha\rightarrow 0$ is equal to
$inf\{x^{*}(y): x^{*}\in S(x)\}$. Therefore, 
the function $\phi (\alpha)=\|x+\alpha y\|$
is differentiable at a point $\alpha \in R$ if and only 
if $x_{1}^{*} (y)=x_{2}^{*} (y)$
for every $x_{1}^{*}, x_{2}^{*}\in S(x+\alpha y).$

Fix linearly independent vectors  $x,y \in E.$ The function 
$\phi (\alpha)=\|x+\alpha y\|$ is convex on $R$ and, hence, $\phi$ is
differentiable almost everywhere on $R$ with respect to Lebesgue measure
(see [5]).

Denote by $D(x,y)$ the set of points $\alpha$ at which $\phi$ is 
differentiable.

\proclaim{Lemma 1} Let $\alpha\in D(x,y), a,b\in R.$ Then

(i)  the number $x^{*}(ax+by)$ does not depend on the choice of 
$x^{*} \in S(x+\alpha y),$

(ii) $x+\alpha y\bot ax+by$ if and only if  $x^{*}(ax+by)=0$ for every
$x^{*} \in S(x+\alpha y).$\endproclaim

\demo{Proof}  (i) As shown above, $x^{*}(y)$ does not depend on the choice of
$x^{*} \in S(x+\alpha y).$ Besides,
$$x^{*}(x) = x^{*}(x+\alpha y) - \alpha x^{*}(y) = \|x+\alpha y\| -
\alpha x^{*}(y)$$
for every $x^{*} \in S(x+\alpha y),$ so $x^{*}(x)$ does not depend on the
choice of a functional $x^{*}$ also.

(ii) If  $x+\alpha y\bot ax+by$ then, by the definition of orthogonality and 
(1), we have 
$\sup\{x^{*}(ax+by): x^{*} \in S(x+\alpha y)\}\geq 0$ and 
$\inf\{x^{*}(ax+by): x^{*} \in S(x+\alpha y)\}\leq 0.$
By (i),  $sup = inf = 0$. 

On the other hand, if $x^{*}(ax+by)=0$ for any  $x^{*} \in S(x+\alpha y)$
then
$$ x^{*}((x+\alpha y) + \gamma (ax+by)) =  x^{*}(x+\alpha y) = 
\|x + \alpha y\|$$
for every $\gamma \in R.$ Since $\|x^{*}\|=1$  we have 
$x+\alpha y\bot ax+by.$\bull\enddemo

The following fact is an easy consequence of the convexity of the function
$\alpha \rightarrow \|x+\alpha y\|$.

\proclaim{Lemma 2} The set of numbers $\alpha$ for which $x+\alpha y \bot y$
is a closed segment $[m,M]$ in $R$  and  $\|x+\alpha y\| = \|x+my\|$ for every
$\alpha \in [m,M]$.\endproclaim

\proclaim{Lemma 3} Let $\alpha\in D(x,y).$ Then either $x+\alpha y \bot y$
or there exists a unique number $f(\alpha)\in R$ such that $x+\alpha y \bot
x - f(\alpha) y.$\endproclaim

\demo{Proof} By Lemma 1, the numbers $x^{*}(x)$ and  $x^{*}(y)$ does not 
depend
on the choice of $x^{*} \in S(x+\alpha y).$ Fix $x^{*} \in S(x+\alpha y).$
If $x^{*}(y)=0$ then, by Lemma 1, $x+\alpha y \bot y$. If $x^{*}(y)\neq 0$
then, again by Lemma 1, $x+\alpha y \bot x -\beta y$ if and only if
$x^{*}(x - \beta y) = 0.$ Thus, $f(\alpha)=x^{*}(x)/x^{*}(y).$\bull\enddemo

By Lemma 2, the function $f$ is defined on $R\setminus [m,M]$. It appears 
that the norm can be expressed in terms of the function $f.$

\proclaim{Lemma 4} For every $\alpha > M,$
$$\|x+\alpha y\| = \|x+My\| \exp(\int_{M}^{\alpha} 
(t + f(t))^{-1} dt) \tag 2$$
and, for every $\alpha < m,$
$$\|x+\alpha y\| = \|x+my\| \exp(-\int_{\alpha}^{m} 
(t + f(t))^{-1} dt) \tag 3$$
\endproclaim

\demo{Proof} Let $\alpha\in D(x,y), \alpha > M.$ Fix  
$x^{*} \in S(x+\alpha y).$
By Lemma 3, $x^{*}(x)=f(\alpha)x^{*}(y)$ and, by (1), 
$x^{*}(y)=\|x+\alpha y\|_{\alpha}^{'}.$
Therefore, $x^{*}(x) = x^{*}(x+\alpha y) - \alpha x^{*}(y) = \|x+\alpha y\| -
\alpha \|x+\alpha y\|_{\alpha}^{'}.$ We have
$$\|x+\alpha y\|_{\alpha}^{'}/\|x+\alpha y\| = (\alpha + f(\alpha))^{-1}$$
Since $\alpha$ is an arbitrary number from $D(x,y)\cap [M,\infty]$ and
Lebesgue measure of the set $R\setminus D(x,y)$ is zero, we get
$$\int_{M}^{\alpha} (\|x+ty\|_{t}^{'}/\|x+ty\|) dt =  
\int_{M}^{\alpha} (t + f(t))^{-1} dt \tag 4$$
for every $\alpha > M.$ It is easy to see that the function
$\alpha \rightarrow ln\|x+\alpha y\|$ satisfies the Lipschitz 
condition and, therefore,
is absolutely continuous. Every absolutely continuous function coincides with
the indefinite integral of its derivative  [5], so the integral in the 
left-hand side of (4) is equal to $ln(\|x+\alpha y\|/\|x+My\|)$ and 
we get (2).
The proof of (3) is similar.\bull\enddemo 

Now we can prove the main result.

\proclaim{Theorem} Let $E$ be a real Banach space and $T:E\rightarrow E$ be
a linear operator preserving orthogonality. Then $T=kU$ where $k\in R$ and
$U$ is an isometry.\endproclaim

\demo{Proof} Assume that $T$ is not the zero operator and fix $x\in E$
such that $Tx\neq 0$. Consider an arbitrary $y\in E$ 
such that $x\neq \alpha y$
for every $\alpha\in R.$ Denote by $I_{1}$ and $I_{2}$ the intervals $[m,M]$
corresponding to the pairs of vectors $(x,y)$ and $(Tx,Ty).$

Since $T$ preserves orthogonality we have $I_{1}\subset  I_{2}.$ Let us prove 
that $I_{1} = I_{2}.$ Assume that $I=I_{2}\setminus I_{1}\neq \emptyset$
and consider a number $\alpha\in I$ such that $\alpha\in D(x,y)\cap D(Tx,Ty).$
Since $\alpha\in I_{2}$ we have $Tx+\alpha Ty \bot Ty.$ By Lemma 3, there 
exists a number $f(\alpha)$ such that $x+\alpha y \bot x - f(\alpha) y$ and,
consequently, $Tx+\alpha Ty \bot Tx - f(\alpha) Ty.$ By Lemma 1, for every
functional $x^{*} \in S(Tx+\alpha Ty),$ we have $x^{*}(Ty)=0$ and 
$x^{*}(Tx-f(\alpha)Ty)=0.$ But then $0=x^{*}(x+\alpha y)=\|x+\alpha y\|$
and we get a contradiction.

Thus, the numbers $m, M$ and, obviously, the function $f(\alpha)$ are the
same for both pairs of vectors $(x,y)$ and  $(Tx,Ty).$

The functions $\|x+\alpha y\|$ and $\|Tx+\alpha Ty\|$ are constant and 
non-zero on $[m, M]$, so there exist $k_{1}, k_{2}\in R$ such that
$\|x+\alpha y\|=k_{1}$ and $\|Tx+\alpha Ty\|=k_{2}$ for every 
$\alpha\in [m, M].$

Using this fact and (2), (3) for both pairs of vectors $(x,y)$ and  $(Tx,Ty)$
we get $\|Tx+\alpha Ty\|=(k_{2}/k_{1}) \|x+\alpha y\|$ for every 
$\alpha \in R.$
First put $\alpha =0$ and then divide the latter equality by 
$\alpha$ and tend 
$\alpha$ to infinity. We get $\|Tx\|/\|x\|=\|Ty\|/\|y\|$ for every non-zero
$x,y\in E$ which completes the proof.\bull\enddemo

\subheading{Acknowledgements.} Part of this work was done when I was visiting 
Memphis State University. I am grateful to Professors J.Jamison and A.Kaminska
for fruitful discussions  and hospitality. I express my gratitude to
Professor J.Arazy for helpful remarks.

\vskip2truecm

\subheading{References}

\item{1.}Dunford, N. and Schwartz, J. {\it Linear operators. Vol.1,}
Interscience, New York 1958.

\item{2.}Giles, J.R. {\it Classes of semi-inner product spaces,}   
Trans. Amer. Math. Soc.  129 (1967), 436-446.

\item{3.}James, R.C. {\it Orthogonality and linear functionals in normed
linear spaces,} Trans. Amer. Math. Soc.  61 (1947), 265-292.

\item{4.}Koehler, D. and Rosenthal, P. {\it On isometries of normed 
linear spaces,}
Studia Mathematica  36 (1970), 213-216.

\item{5.}Riesz, F. and Sz-Nagy, B. {\it Lecons d'analyse fonctionelle,}
Akademiai Kiado, Budapest 1972.
\bye